\documentclass[review]{elsarticle}

\usepackage{lineno,hyperref}

\usepackage{epsfig}
\usepackage{amsmath}
\usepackage{amsfonts}
\usepackage{amsthm}
\usepackage{amsbsy}
\usepackage{appendix}
\usepackage{epstopdf}
\usepackage[table]{xcolor}
\usepackage{subfloat}
\usepackage{float}
\usepackage[thinlines]{easytable}
\usepackage{array}
\usepackage{color}
\usepackage{inputenc}

\numberwithin{equation}{section}
\allowdisplaybreaks

\theoremstyle{plain}

\newtheorem{algo}{Algorithm}[section]
\newtheorem{assump}{\bf Assumption}[section]

\newtheorem{definition}{\bf Definition}[section]
\theoremstyle{remark}
\newtheorem{remark}{\bf Remark}[section]
\theoremstyle{remark}

\modulolinenumbers[5]

\makeatletter
\def\ps@pprintTitle{%
   \let\@oddhead\@empty
   \let\@evenhead\@empty
   \let\@oddfoot\@empty
   \let\@evenfoot\@oddfoot
}
\makeatother

\begin{document}

 \graphicspath{
               {Figures/}
              }

\begin{frontmatter}

\title{Dynamic Mode Decomposition for Construction of Reduced-Order Models of Hyperbolic Problems with Shocks }
\author[mymainaddress]{Hannah Lu\fnref{email1}}
\fntext[email1]{email: \texttt{hannahlu{@}stanford.edu}}
\author[mymainaddress]{Daniel M. Tartakovsky\corref{mycorrespondingauthor}}
\cortext[mycorrespondingauthor]{Corresponding author}
\ead{tartakovsky@stanford.edu}
\address[mymainaddress]{Department of Energy Resources Engineering, Stanford University, Stanford, CA 94305, USA}

\begin{abstract}
Construction of reduced-order models (ROMs) for hyperbolic conservation laws is notoriously challenging mainly due to the translational property and  nonlinearity of the governing equations. While the Lagrangian framework for ROM construction resolves the translational issue, it is valid only before a shock forms. Once that occurs, characteristic lines cross each other and projection from a high-fidelity model space onto a ROM space distorts a moving grid, resulting in numerical instabilities. We address this grid distortion issue by developing a  physics-aware dynamic mode decomposition (DMD) method based on hodograph transformation. The latter provides a map between the original nonlinear system and its linear counterpart, which coincides with the Koopman operator. This strategy is consistent with the spirit of physics-aware DMDs in that it retains information about shock dynamics. Several numerical examples are presented to validate the proposed physics-aware DMD approach to constructing accurate ROMs.
\end{abstract}

\begin{keyword}
Conservation law \sep  Koopman operator \sep Proper orthogonal decomposition \sep Hodograph transformation
\end{keyword}

\end{frontmatter}

\linenumbers

\section{Introduction}

Since introduction of Euler equations, hyperbolic conservation laws play a significant role in gas dynamics, astrophysics, plasma, traffic flow, multiphase flow in porous media~\cite{chang1989riemann,shu1991physics,whitham2011linear,courant1999supersonic,bear2013dynamics} and other fields of science and engineering. Wave-like solutions of hyperbolic equations can exhibit various rarefaction and shock behaviors, whose occurrence  strongly depends on a functional form of the flux function. Discontinuity and uniqueness of such solutions pose challenges in theoretical treatment of hyperbolic conservation laws~\cite{oleinik1957discontinuous,smoller2012shock}. Theoretical advances, such as entropy conditions and the concept of a weak solution~\cite{lax1971shock,harten1983symmetric}, ameliorate this difficulty by providing physical interpretation to these solutions. Likewise, numerical  high-resolution methods have been designed to resolve nonlinearities and accurately capture shocks~\cite{leveque1992numerical, osher1984high, majda1979numerical}.  Although continued developments in scientific computing have improved the performance of high-resolution simulations, their computational cost is often too high to model complex systems at spatiotemporal resolutions and scales of interest. The cost can become prohibitive when used in the context of uncertainty quantification or data assimilation, both of which require a large number of repeated forward model runs.

Reduced-order models (ROMs) provide an efficient alternative to their high-fidelity, physics-based counterparts that can be deployed in large-scale multiphysics simulations. Robust tools for construction of ROMs for problems described by ordinary differential equations or parabolic partial differential equations (PDEs) include proper orthogonal decomposition (POD)~\cite{benner2015survey, kerschen2005method, rowley2005model} and dynamic mode decomposition (DMD)~\cite{schmid2010dynamic, alla2017nonlinear, brunton2016koopman, williams2013hybrid}. 
The challenge of extending these techniques to hyperbolic or advection-dominated parabolic PDEs with smooth solutions was met in~\cite{lu2019part1} through development of the physics-aware DMD and POD approaches within a Lagrangian framework.  However, in the presence of strong shocks and/or sharp gradients, Lagrangian POD methods can generate numerical instability caused by grid distortion~\cite{mojgani2017lagrangian}. Once characteristics of a nonlinear hyperbolic PDE intersect each other, the projection from a high-dimensional manifold of the high-fidelity model (HFM) onto a low-dimensional subspace of the low-fidelity model (LFM), e.g., ROM, is not guaranteed and typically fails to preserve topological properties of the original HFM. We elaborate on this point in section~\ref{sec:ROMfailure}, in terms relevant to DMD.

We use hodograph transformation~\cite{hamilton1847hodograph} to resolve this outstanding issue in construction of ROMs for PDEs with discontinuous solutions and shocks. Hodograph diagrams have originated in meteorology to plot wind from soundings of the Earth's atmosphere. Since then,  hodograph transformation morphed into a technique designed to transform nonlinear PDEs into linear ones by  interchanging the dependent and independent variables in a PDE. Hodograph-type transformations have been used to find quasilinear analogues of  semi-linear equations, and to derive new analytical solutions to special classes of PDEs~\cite{clarkson1989hodograph}. 
Advantages of mapping nonlinear PDEs onto their linear counterparts are self-evident: analytical tools available for linear PDEs provide better understanding of a solution's behavior, and numerical solvers for linear systems are both easier to implement and computationally cheaper. 

The Koopman operator theory~\cite{koopman1931hamiltonian} shares the goal of hodograph transformation: a Koopman operator is an infinite-dimensional linear operator that represents the underlying finite-dimensional nonlinear dynamic system by judiciously choosing observable functions. It is also similar in its goal to integral transformations that map certain classes of nonlinear PDEs onto their linear counterparts. For example, the Cole-Hopf transformation and the Kirchhoff transformation map, respectively, Burgers' equation and a class of nonlinear diffusion (heat conduction) equations onto a linear diffusion equation. 
These integral transformations have been used in the context of the Koopman operator theory and DMD/POD to constructed ROMs for Burgers' equation~\cite{nathan2018applied} and a nonlinear diffusion equation~\cite{lu2019error}. A major goal of our study is to establish clear connections between hodograph transformation and the Koopman theory. This relationship between the two is then used both to identify observables for a Koopman operator via hodograph transformation and to construct ROMs for hyperbolic conservation laws with shocks via DMD.

A general procedure of the physics-aware DMD algorithm and its connection with Koopman operators are reviewed in section~\ref{sec:ROMfailure}. This section also contains a numerical demonstration of the failure of Lagrangian-based ROMs to capture the dynamics described by conservation laws with shocks. We illustrate the use of hodograph transformation by analyzing the inviscid Burgers' equation (\ref{sec:hodograph}) and more general hyperbolic PDEs with a convex flux function (appendix~\ref{app:Buckley}). In section~\ref{sec:DMD}, we combine hodograph transformation with the Koopman operator theory to design a physics-aware DMD algorithm for construction of ROMs for conservation laws with shocks. Several numerical tests are presented in section~\ref{sec:examples} to validate the proposed physics-aware DMD approach. Main conclusions drawn from our study are summarized in section~\ref{sec:conc}.

\section{Construction of ROMs and their Failure for Problems with Shocks}
\label{sec:ROMfailure}

Consider a state variable $u(x,t): [a,b] \times [0,T] \rightarrow \mathbb R$, where the constants $a, b, T \in \mathbb R$. The dynamics of $u(x,t)$ is described by a one-dimensional scalar conservation law 
\begin{equation}\label{2-2}
\frac{\partial u}{\partial t}+\frac{\partial F(x,t,u)}{\partial x} =0
\quad\text{or}\quad 
\frac{\partial u}{\partial t}+f(x,t,u)\frac{\partial u}{\partial x} =0, \quad f(x,t,u) = \frac{\partial F(x,t,u)}{\partial u}.
\end{equation}
This hyperbolic PDE is subject to the initial condition $u(x,t=0) = u_0(x)$ and, when appropriate (i.e., when $|a|, |b| < \infty$), boundary conditions at $a$ and/or $b$. 
%
%
The intervals $[0,T]$ and $[a,b]$ are discretized with $(N+1)$ and $J$ nodes separated, respectively, by $\Delta t$ and $\Delta x$. To be specific, we solve~\eqref{2-2} with a conservative first-order upwind scheme~\cite{leveque1992numerical} 
\begin{equation}\label{2-3}
u_j^{n+1} = u_j^n-\frac{\Delta t}{\Delta x}(F_{j+1/2}^n-F_{j-1/2}^n),
\end{equation}
where $n = 0,\ldots,N$ indicates the $n$th time step, with $n=0$ corresponding to $t=0$, and $n=N$ to $t = T$;  $j = 1,\ldots,J$ denotes the $j$th spatial node, such that $j=1$ and $J$ coincide with $x = a$ and $b$, respectively; and
$$\begin{aligned}
&F_{j+1/2}^n = \frac{ F(\cdot, u_{j+1}^n)+F(\cdot, u_j^n) }{2} - |a_{j+1/2}^n| \frac{u_{j+1}^n-u_j^n}{2},\\
&a_{j+1/2}^n = \left\{
\begin{aligned}
&\frac{F_{j+1}^n-F_j^n}{u_{j+1}^n-u_j^n}&&\mbox{if}&u_{j+1}^n\neq u_j^n,\\
&f(\cdot, u_j)&&\mbox{if}&u_{j+1}^n=u_j^n.
\end{aligned}
\right.
\end{aligned}$$
A numerical solution provided by~\eqref{2-3} with sufficiently small $\Delta t$ and $\Delta x$ are referred to as a reference HFM throughout the paper.

Standard (Eulerian) approaches to construction of a ROM for~\eqref{2-2} often fail due to the traveling-wave nature of its solution~\cite{lu2019part1, mojgani2017lagrangian}. In a shock-free scenario, the Lagrangian framework can resolve the translational issue in the POD or DMD approaches to ROMs by keeping track of the characteristic lines. 

In the Lagrangian framework,~\eqref{2-2} becomes
\begin{equation}\label{2-4}
\left\{
\begin{aligned}
&\frac{\text dx}{\text dt} = f(x,t,u), \qquad x(0) = \eta \\
&\frac{\text du}{\text dt} =0, \qquad u(\eta,0) = u_0(\eta),
\end{aligned}
\right.
\end{equation}
where $\eta \in \mathbb R$ is a label of the characteristic $x(t)$. As in the Eulerian case, we use the uniform discretization of the time interval $[0,T]$, such that $0=t^0<t^1<\ldots<t^N=T$ with time step $\Delta t = t^{n+1}-t^n$.  At time $t=0$, the space, $[a,b]$, is discretized with a uniform mesh $\mathbf x^0 = [x^0_1,\ldots, x^0_J]^\top$ of mesh size $\Delta x^0  = x^0_{j+1}-x^0_j$. Unlike its Eulerian counterpart, the spatiotemporal discretization of $u(x,t)$ in the Lagrangian framework, $\mathbf u^n = [u_1^n,\ldots, u_J^n]^\top$ for $n = 0,\ldots,N$, may be nonuniform in space due to the temporal evolution of the grid nodes $x_j(t)$. The backward Euler discretization, used in~\cite{mojgani2017lagrangian}, transforms~\eqref{2-4} into
\begin{equation}\label{2-5}
\left\{
\begin{aligned}
&x_j^{n+1} = x_j^n+\Delta t f(x_j^{n+1}, (n+1) \Delta t, u_j^{n+1}),\\
&u_j^{n+1} = u_j^n
\end{aligned}
\right.
\end{equation}
or, in vector form,
\begin{equation}\label{2-6}
\left\{
\begin{aligned}
&\mathbf R_x(\mathbf x^{n+1}) \equiv \mathbf x^{n+1}-\mathbf x^n-\Delta t\mathbf f^{n+1}(\cdot, \mathbf u ^{n+1}) = \mathbf 0, \\
&\mathbf R_u(\mathbf u^{n+1}) \equiv \mathbf u^{n+1}-\mathbf u^n = \mathbf 0,
\end{aligned}
\right.
\end{equation}
where $\mathbf x^n = [x_1^n, \ldots, x_J^n]^\top$ denotes the nodes of the Lagrangian computational grid at the $n$th time step. This numerical scheme involves $N$ iterations in the two high-dimensional $J \times 1$ vectors, $\mathbf x^{n+1}$ and $\mathbf u^{n+1}$. It provides a Lagrangian HFM. 

To construct a ROM, a data set consisting of a sequence of $M$ solution snapshots ($M \le N$ and, ideally, $M \ll N$) is collected from the HFM. Since $\mathbf u^n$ is conservative and invariant in time, we only need the data matrix $\mathbf X$ with $M$ snapshots of the Lagrangian grid $\mathbf x^n$:
\begin{equation}\label{2-7}
\mathbf X = \begin{bmatrix}
|&|&&|\\
\mathbf x^1&\mathbf x^2&\cdots&\mathbf x^M \\
|&|&&|
\end{bmatrix}.
\end{equation}
In the next two subsections, we briefly revisit the algorithms of Lagrangian POD~\cite{mojgani2017lagrangian} and Lagrangian DMD~\cite{lu2019part1} used to construct a ROM.

\subsection{Lagrangian POD}

Identification of the POD modes is based on a reduced singular value decomposition (SVD),
\begin{equation}\label{2-8}
\mathbf X = \mathbf U\boldsymbol \Sigma \mathbf V^*,
\end{equation}
where $\mathbf U \in \mathbb C ^{J \times K},\boldsymbol \Sigma = \mathbb C ^{K\times K},\mathbf V \in \mathbb C^{M \times K}$, $K$ is the rank of the matrix $\mathbf X$ approximated by the reduced SVD. Further rank truncation can be achieved by using the energy criterion,
\begin{equation}\label{2-9}
r = \min_{k}\left\{\frac{\sigma_k}{\sum_{k=1}^K \sigma_k}<\varepsilon\right\},
\end{equation}
where $\sigma_k$ are the diagonal elements of $\boldsymbol \Sigma$, and $\varepsilon$ is a small number (tolerance), chosen to be $\varepsilon = 10^{-4}$ in all our numerical examples. After the truncation, one gets the POD modes
\begin{equation}\label{2-10}
\boldsymbol \Phi = \mathbf U(:,1:r) =  \begin{bmatrix}
|&|&&|\\
\phi_1&\phi_2&\cdots&\phi_r\\
|&|&&|
\end{bmatrix}.
\end{equation} 
Notice that $r\ll K\leq\min\{J,M\}$, and the basis $\{\phi_1,\cdots,\phi_r\}$ is orthonormal. Galerkin projection in the low-dimensional space spanned by the POD basis provides a ROM (low-fidelity solution),
\begin{equation}\label{2-11}
\mathbf x_\text{POD}^{n+1} = \sum_{k = 1}^r\hat x_k^{n+1}\phi_k = \boldsymbol \Phi\hat{\mathbf x}^{n+1}.
\end{equation}
The $r\times 1$ vector $\hat{\mathbf x}^{n+1}$ of coefficients is computed as a solution of 
\begin{equation}\label{2-12}
\boldsymbol \Phi^\top\mathbf R\left(
\boldsymbol \Phi
\begin{bmatrix}
|\\
\mathbf {\hat x}^{n+1}\\
|
\end{bmatrix}
\right)=0.
\end{equation}
that is obtained by substituting~\eqref{2-11} into the first equation in~\eqref{2-6} and projecting onto the subspace spanned by $\boldsymbol \Phi$.

\subsection{Lagrangian DMD}

The DMD algorithm~\cite{lu2019part1} is applied to the Lagrangian grid matrix $\mathbf X$ in~\eqref{2-7}.
\begin{algo}{Lagrangian DMD algorithm}
\begin{itemize}
\item[0.] Create data matrices of $(M-1)$ observables, $\mathbf X_1$ and $\mathbf X_2$,
\begin{equation}\label{2-13}
\mathbf X_1 = \begin{bmatrix}
|&|&&|\\
\mathbf x^1&\mathbf x^2&\cdots&\mathbf x^{M-1}\\
|&|&&|
\end{bmatrix}, \qquad 
\mathbf X_2 = \begin{bmatrix}
|&|&&|\\
\mathbf x^2&\mathbf x^3&\cdots&\mathbf x^{M}\\
|&|&&|
\end{bmatrix}.
\end{equation}
\item[1.] Apply SVD of matrix $\mathbf X_1 \approx \mathbf U\boldsymbol \Sigma \mathbf V^*$ with $\mathbf U \in \mathbb C^{J \times r}, \boldsymbol\Sigma \in \mathbb C^{r\times r}, \mathbf V\in \mathbb C^{r\times M}$, where $r$ is the truncated rank chosen by a certain criterion, e.g.,~\eqref{2-9}.
\item [2.] Compute $\tilde{\mathbf K}=\mathbf U^*\mathbf X'\mathbf V\boldsymbol\Sigma^{-1}$ as an $r\times r$ low-rank approximation of $\mathbf K$.
\item [3.] Compute eigen-decomposition of $\tilde{\mathbf K}$: $\tilde{\mathbf K} \mathbf W = \mathbf W\boldsymbol\Lambda$, $\boldsymbol\Lambda = (\lambda_k)$.
\item [4.] Reconstruct eigen-decomposition of $\mathbf K$. Eigenvalues are $\boldsymbol\Lambda$ and eigenvectors are $\boldsymbol\Phi  = \mathbf U\mathbf W$.
\item [5.] Future $\mathbf x_\mathrm{DMD}^{n+1}$ is predicted by
\begin{equation}\label{2-14}
\mathbf x_\mathrm{DMD}^{n+1} =\boldsymbol \Phi\boldsymbol\Lambda^{n+1}\mathbf b, \ \ n > M
\end{equation}
with $\mathbf b =\boldsymbol \Phi^{-1}\mathbf x_1$.
\item [6.] Interpret the solution in the moving grid:
\begin{equation}\label{2-15}
 u_\mathrm{DMD}(x_j^n, n \Delta t) = u_0(x_j^0).
\end{equation}
\end{itemize}
\end{algo}

\subsection{ROM Failure for Problems with Shocks: Inviscid Burgers Equation}
\label{sec:DMDfailure}

One of the most studied examples of~\eqref{2-2} is the inviscid Burgers equation
\begin{equation}\label{2-16}
\frac{\partial u}{\partial t}+u\frac{\partial u}{\partial x} = 0, \qquad u(x,0)=u_0(x);
\end{equation}
which we define on the space-time domain $(x,t) \in [0,2\pi]\times [0,1]$. Depending on the boundary and initial conditions, this problem admits both smooth and discontinuous solutions $u(x,t)$. For example, a smooth solution is obtained for the periodic boundary conditions, $u(0,t) = u(2\pi,t)$, and the initial data $u_0(x) = 1+\sin(x)$. In this setting, standard (Eulerian) ROMs fail, while the ROMs based on either Lagrangian POD or Lagrangian DMD perform well in terms of both accuracy and computational efficiency~\cite{lu2019part1}.


%



\begin{figure}[tphb]
\includegraphics[width=0.5\textwidth]{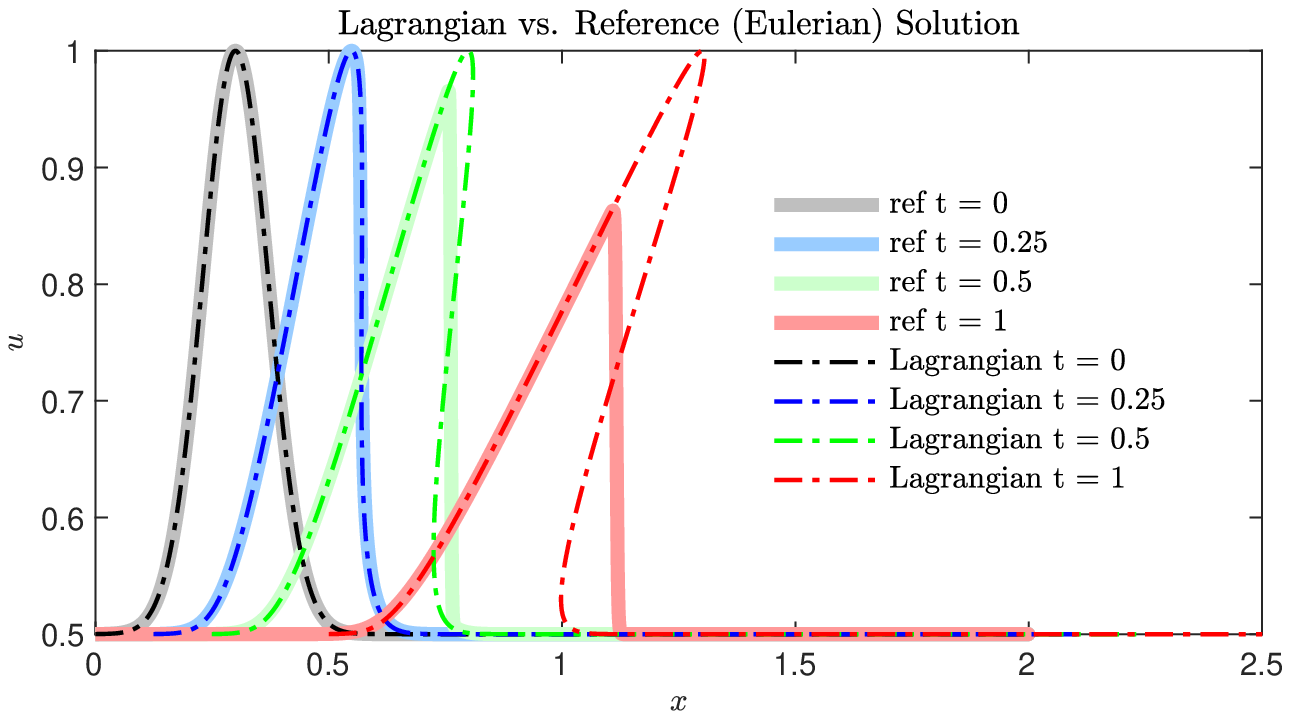}
\includegraphics[width=0.5\textwidth]{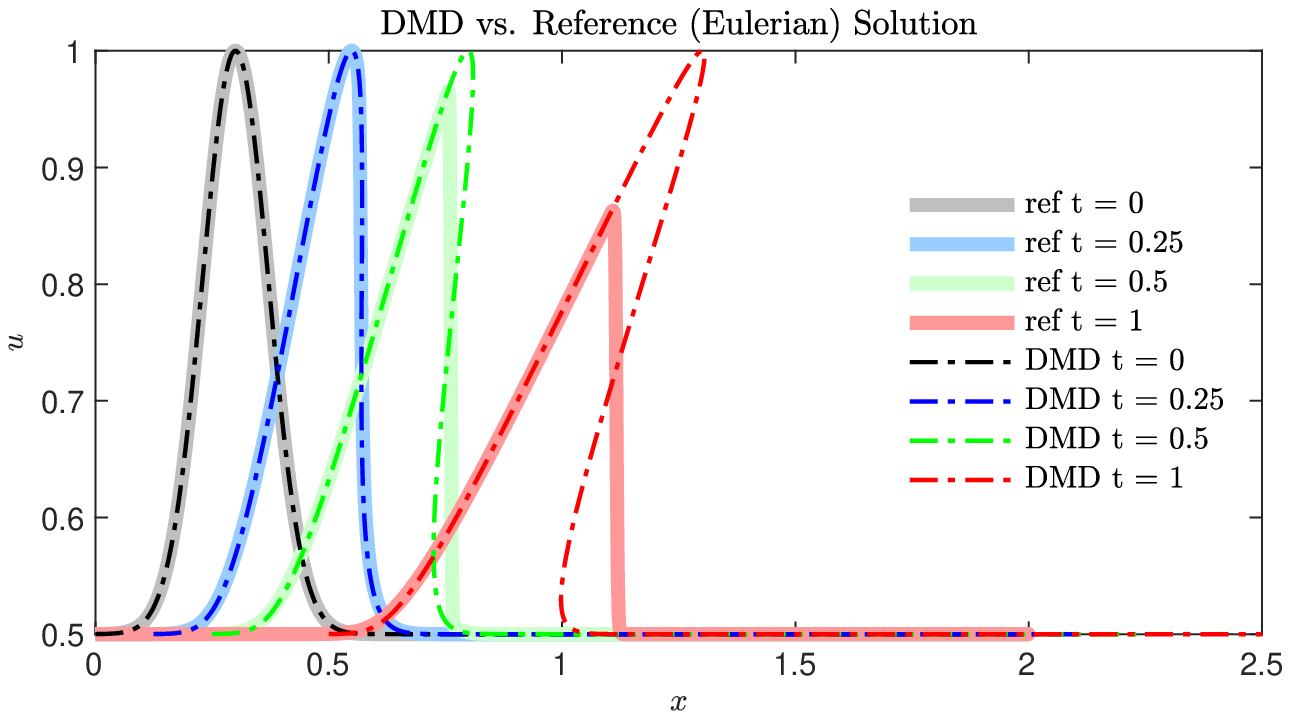}
\caption{Solutions of the inviscid Burgers equation with a shock. \textbf{(a)} The full solution obtained with the Lagrangian numerical scheme~\eqref{2-18} leading to the overshoot. \textbf{(b)} The Lagrangian DMD solution trained on a few snapshots of the faulty full solution. The reference solution is obtained with~\eqref{2-3}.}
\label{fig:badLagrange}
\end{figure}

A solution to~\eqref{2-16} develops shocks in finite time for, e.g., a Gaussian-type initial data, 
\begin{equation}\label{2-17}
u_0(x) = \frac{1}{2} + \frac{1}{2} \exp \! \left[- \frac{(x-0.3)^2}{0.01} \right].
\end{equation}
In the pure Lagrangian approach~\eqref{2-4}, the discretization has to account for shock formation. Once the characteristic lines cross each other, the Lagrangian mesh becomes sensitive to the choice of discretization of $u[x(t),t]$. 
For instance, a discretization of~\eqref{2-4} with $f(\cdot, u) =u$,
\begin{equation}\label{2-18}
\left\{
\begin{aligned}
&u_j^{n+1} = u_j^n,\\
&x_j^{n+1} = x_j^n +\frac{\Delta t}{2}(u_j^n+u_j^{n+1}),
\end{aligned}
\right.
\end{equation}
would lead to the so-called ``overshoot" that admits multi-value solutions (Figure~\ref{fig:badLagrange}(a)), which contradicts the entropy condition. This is a typical problem with the Lagrangian framework. It should come as no surprise that an attempt to build a ROM with the Lagrangian DMD based on the faulty discretization~\eqref{2-18} likewise results in failure (Figure~\ref{fig:badLagrange}(b)). The Lagrangian DMD faithfully reproduces the unphysical solution obtained with the faulty discretization scheme~\eqref{2-18}. In other words, the resulting unphysical ROM is not caused by the DMD algorithm itself; the data from the full Lagrangian model~\eqref{2-18} provide inaccurate and incomplete (without shock) information from the very beginning.

\begin{figure}[tphb]
\includegraphics[width=0.5\textwidth]{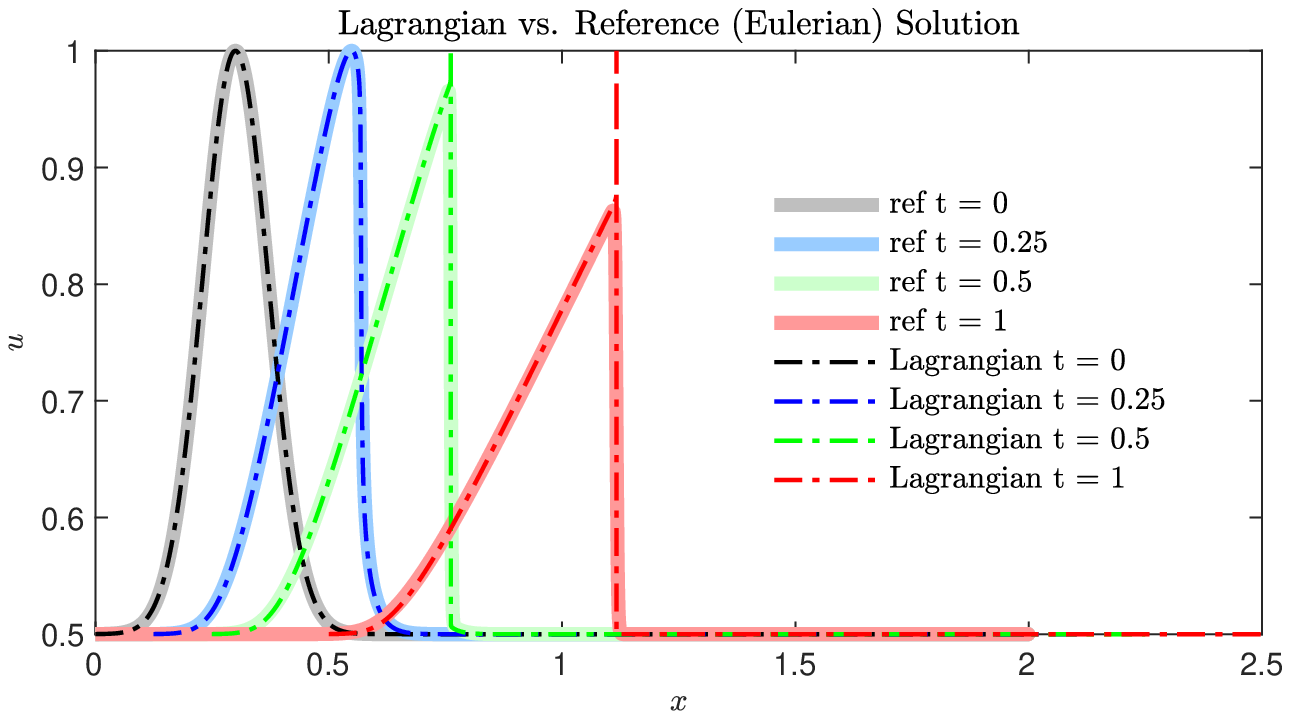}
\includegraphics[width=0.5\textwidth]{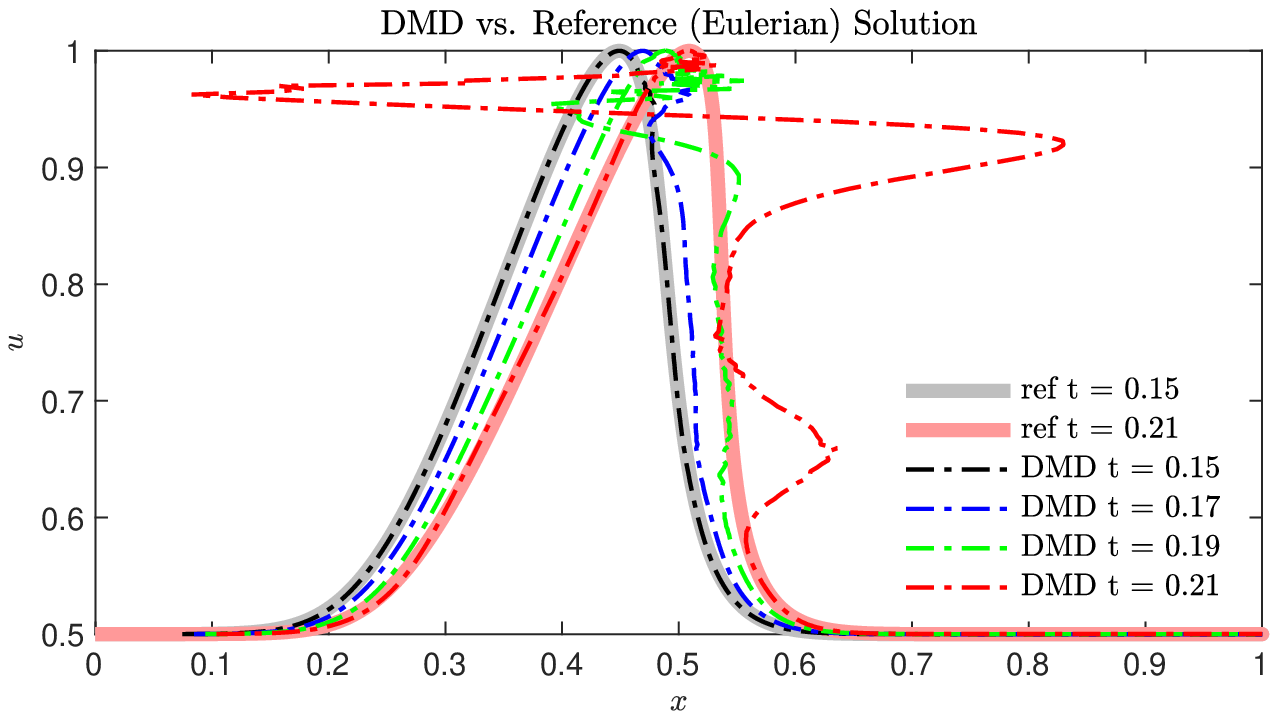}
\caption{Solutions of the inviscid Burgers equation with a shock. \textbf{(a)} The full solution obtained with the appropriate Lagrangian numerical scheme~\eqref{2-19}. \textbf{(b)} The Lagrangian DMD solution trained on a few snapshots of the accurate full solution. The reference solution is obtained with~\eqref{2-3}.}
\label{fig:goodLagrange}
\end{figure}

To isolate the performance of the Lagrangian DMD, we consider a numerical scheme that is known for its ability to handle shocks: the backward semi-Lagrangian method (BSLM)
\begin{equation}\label{2-19}
\left\{
\begin{aligned}
&u_j^{n+1} = u_j^n,\\
&x^* = x_j^n+\frac{\Delta t}{2}u_j^n,\\
&x_j^{n+1} = x_j^n +\frac{\Delta t}{2}(u(x^*,t^n)+u(x^*,t^{n+1})).
\end{aligned}
\right.
\end{equation}
We employed the (explicit) mid-point rule to avoid implicit iterations. As one can see from the above figure, (\ref{2-19}) gives the right physical solution with shock features. Figure~\ref{fig:goodLagrange}(a) reveals that this numerical scheme is indeed capable to accurately approximate the solution of the inviscid Burgers equation with shocks. However, the Lagrangian DMD algorithm using snapshots from the full solution~\eqref{2-19} suffers from instability once a shock is about to form (Figure~\ref{fig:goodLagrange}(b)). The grid becomes severely distorted once the characteristic lines intersect each other at the interface where sharp gradients of $u(x,t)$ occur.  At the intersect, one arrival location of $x$ corresponds to two different departure values of $u$. However, the DMD modes projection from the HFM to the ROM does not keep the topological information about this multivalued mapping in the ROM process, resulting in the Lagrangian grid distortion. 

\begin{remark}
The Lagrangian POD approach suffers from similar problems~\cite{mojgani2017lagrangian}. Moreover, the POD projection on the accurate Lagrangian scheme (\ref{2-19}) would still require interpolation in the high-dimensional space. One might need techniques such as DEIM~\cite{chaturantabut2010nonlinear} to keep the resulting ROM's efficiency. Nevertheless, extensions of POD are beyond the scope of our study; we focus on DMD-based ROMs due to their iteration-free nature.
\end{remark}

\section{Hodograph Transformation}
\label{sec:hodograph}

We start with a mathematical definition of hodograph transformation reproduced from~\cite{clarkson1989hodograph}.
\begin{definition}\label{def:hodograph}
A pure hodograph transform is a transformation of the form
\begin{equation}\label{3-1}
\tau = t, \quad \xi = u(x,t).
\end{equation}
\end{definition}
%
%
For the inviscid Burgers equation~\eqref{2-16}, we first consider a scenario where only one shock is developed from the initial data $u_0(x)$ in finite time. This necessitates the following assumption.
\begin{assump}\label{assumption1}
The function $u_0(x)$ satisfies four conditions:
\begin{itemize}
\item $u_0(x)$ is smooth;
\item $u_0(x)$ decreases monotonically, i.e., $u_0^\prime(x) <0$ for all $x$; and 
 $\lim_{x\to +\infty} u_0(x) = u_R$, $\lim_{x\to -\infty}u_0(x) = u_L$ with constants $u_R < u_L$;
\item $u_0(x)$ has a unique inflection point $(x^*, u^*)$ with $u^* = u_0(x^*)$, meaning $u_0^{\prime\prime}(x^*) = 0$;
\item $u_0^{\prime\prime\prime} (x^*) >0$.
\end{itemize}
\end{assump}
This assumption ensures existence of an inverse function $x(t,u) : [0,T]\times [u_R,u_L] \to [a,b]$ of the monotonic function $u(t,x) : [0,T] \times [a,b] \to [u_R,u_L] $. 
It follows from Definition~\ref{def:hodograph} that the inverse function $x(t,u)$ is a pure hodograph transform. 

\subsection{Solution Before Shock Formation}
\label{sec:3-1}

With $u$ acting as the independent variable and $x$ as the dependent variable, hodograph transformation $x = x(t,u)$ maps the equation for characteristics~\eqref{2-4} of the inviscid Burgers equation~\eqref{2-16}, before the shock formation time $t^*$ (defined later), onto  
\begin{equation}\label{3-5}
\frac{\text d x}{\text dt}(t,u) = u, \quad x(0,u) = x_0(u); \qquad\text{for}\quad (t,u) \in [0,t^*) \times [u_R,u_L].
\end{equation}
Assumption~\ref{assumption1} translates into conditions on the function $x_0(u)$:
\begin{itemize}
\item $x'_0 (u)<0$;
\item $x_0(u)$ has a unique inflection point at $(u^*,x^*)$;
\item $x'''_0(u^*)<0$.
\end{itemize}
Differentiation of~\eqref{3-5} with respect to $u$ gives
\begin{equation}\label{3-6}
\frac{\partial^2 x}{\partial t \partial u}(t,u) = 1,
\end{equation}
from which
\begin{equation}\label{3-7}
\frac{\partial x}{ \partial u}(t,u) = x'_0(u) + t, \qquad\text{for}\quad (t,u) \in [0,t^*) \times [u_R,u_L].
\end{equation}
Since $x'_0(u)<0$, 
$\partial_u x (t,u) < 0$ as long as 
$t^* = \min_u[-x'_0(u)] = -x'_0(u^*)$.
Thus defined $t^*$ determines the time of shock formation. The shock location is $x^* = x(u^*,x^*)$.

\subsection{Solution After Shock Formation}

At times $t$ larger than $t^* = -x'_0(u^*)$, i.e., once the shock forms,~\eqref{3-5} is no longer valid. In the $(x,u)$ plane, one would use the entropy (Rankine-Hugoniot) condition to construct a weak formulation of Burgers' equation. Its analog in the $(u,x)$ plane gives an equation for the shock speed $s$: 
\begin{equation}\label{3-10}
s = \frac{1}{2} \frac{u_1^2 - u_2^2}{u_1-u_2} =\frac{u_1+u_2}{2},
\end{equation}
where $u_1(t)$ and $u_2(t)$ are defined as the limits of $u(t)$ from the top and bottom of the shock, respectively. They are computed as solutions of a system of ordinary differential equations~\cite{ruiwen2018}
\begin{equation}\label{3-12}
\left\{
\begin{aligned}
&\frac{\text du_1}{\text dt} =\frac{1}{2}\frac{u_1-u_2}{g(u_1)-t},\\
&\frac{\text du_2}{\text dt} =-\frac{1}{2}\frac{u_1-u_2}{g(u_2)-t},
\end{aligned}
\right.
\end{equation}
where $g(u) \equiv -x'_0(u)$. These ODEs are subject to initial conditions $u_1(t^*) =u^*$ and $u_2(t^*) = u^*$.
Since $s = \text d x^* / \text dt$, an equation for the shock trajectory $x^*(t)$ is
\begin{equation}\label{3-11}
\frac{\text d x^* }{\text dt} =\frac{u_1+u_2}{2}.
\end{equation}

\subsection{Summary of Hodograph Solution}

Under Assumption~\ref{assumption1}, the hodograph-transformed Burgers equation~\eqref{2-16} in the Lagrangian framework~\eqref{2-4} takes the form of the following ODEs for $x(t,u)$: 
\begin{equation}\label{3-13}
\left\{
\begin{aligned}
&t<t^*: &&\mbox{Equation~\eqref{3-5}}\\
&t>t^*:&&\left\{
\begin{aligned}
&\mbox{Equation (\ref{3-5})}&&\mbox{for}&&u\in(u_R,u_2)\cup(u_1,u_L)\\
&\mbox{Equation (\ref{3-11})}&&\mbox{for}&&u\in(u_2,u_1),
\end{aligned}
\right.
\end{aligned}
\right.
\end{equation}
where $t^* = -x'_0(u^*)$, and $u_1$ and $u_2$ are solutions of~\eqref{3-12}.

\begin{remark}
One can show that $u_1(t)$ is a monotonically increasing function and $u_2(t)$ is a monotonically decreasing function, such that
\begin{equation}\label{3-14}
u_1\geq u^*, \quad u_2\leq u^*, \quad x'_0(u_1)+t\leq 0, \quad x'_0(u_2)+t\leq 0.
\end{equation}
In many cases of interests, and in all our numerical experiments, either $u_2 = u_R$ and $u_1 = u_L$ or $|u_2 - u_1| \ll \Delta t$ (so that $u_2 \approx u_R$ and $u_1 \approx u_L$). This allows one to focus on shock propagation, i.e., on~\eqref{3-11}, without having to solve~\eqref{3-12}.
\end{remark}


\begin{remark}
Functions $u_0(x)$ that do not satisfy Assumption~\ref{assumption1}, such as~\eqref{2-17}, require a decomposition of the initial data into monotonic parts. Each monotonic piece of $u_0(x)$ would have a unique inverse function $x_0(u)$. The entropy condition implies that the increasing $x_0(u)$, i.e. $x'_0(u)>0$, results in a rarefaction solution, which satisfies~\eqref{3-5}. The union of  the rarefaction pieces and shock pieces would give the full solution.
\end{remark}

\begin{remark}
The inviscid Burgers equation is an example of hyperbolic conservation laws with monotonically increasing flux functions $F(\cdot, u)$. Generalization to hyperbolic conservation laws with a convex flux is presented in appendix~\ref{app:Buckley}.
\end{remark}

\section{Physics-Aware DMD for Conservation Laws with Shocks}
\label{sec:DMD}

Previous theoretic investigations, e.g.,~\cite{lu2019error}, demonstrated that the key to a DMD's success in capturing nonlinear dynamics is to identify the underlying Koopman operator. Several numerical studies~\cite{nathan2018applied, rowley2009spectral, williams2015data} confirmed this finding. The Koopman operator theory ensures that a DMD algorithm utilizes all relevant physical information to learn the dynamics. We refer to this approach as physics-aware DMD to distinguish it from the conventional DMD that learns only from (simulations-generated) data.

We review the Koopman operator theory and analyze its connection with hodograph transformation in the context of hyperbolic conservation laws.  Then, we present our general framework for physics-aware DMD for problems with discontinuous solutions and shocks.

A suitable spatial discretization of~\eqref{2-2} leads to a nonlinear dynamical system 
\begin{equation}\label{4-1}
\frac{ \text d\mathbf u}{ \text dt} = \mathcal N(\mathbf u),
\end{equation}
where $\mathbf u(t) = (u_1,\ldots,u_J)^\top \in \mathcal M \subset \mathbb R^J$ is the solution vector with $u_j(t) = u(x_j,t)$ and discretization nodes $x_j$ ($j=1,\ldots,J$); and $\mathcal N$ is a finite-dimensional nonlinear operator.  A flow map $\mathcal N_t :\mathcal M \to \mathcal M$,
\begin{equation}\label{4-2}
\mathcal N_t [\mathbf u(t_0)] \equiv \mathbf u(t_0+t) = \mathbf u(t_0)+\int_{t_0}^{t_0+t} \mathcal N[\mathbf u(\tau)] \text d \tau,
\end{equation}
induces the corresponding discrete-time dynamical system
\begin{equation}\label{4-3}
\mathbf u^{n+1} = \mathcal N_t(\mathbf u^n).
\end{equation}

\begin{definition}[Koopman operator: \cite{kutz2016dynamic}]
Consider a state $\mathbf u$ on a smooth $J$-dimensional manifold $\mathcal M$, whose dynamics are described by~\eqref{4-1}. The Koopman operator $\mathcal K$ is an infinite-dimensional linear operator that acts on all observable functions $g: \mathcal M\to \mathbb C$ such that
\begin{equation}\label{4-4}
\mathcal K g(\mathbf u) = g[\mathcal N(\mathbf u)].
\end{equation}
The discrete-time Koopman operator $\mathcal K_t$ for the discrete dynamical system~\eqref{4-3} is defined as
\begin{equation}\label{4-5}
\mathcal K_t g(\mathbf u^{n}) = g[\mathcal N_t(\mathbf u^n)] =g(\mathbf u^{n+1}).
\end{equation}
\end{definition}
A major benefit of the Koopman operator is the transformation of the finite-dimensional nonlinear problem~\eqref{4-3} in the state space into an infinite-dimensional linear problem (\ref{4-5}) in the observable space. Since $\mathcal K_t$ is an infinite-dimensional linear operator, it has an infinite number of eigenvalues $\{\lambda_k\}_{k=1}^{\infty}$ and eigenfunctions $\{\phi_k\}_{k=1}^\infty$.  In practice, one has to make a finite approximation of the eigenvalues and eigenfunctions. The following assumption is essential to the finite approximation and the choice of observables:

\begin{assump}
Consider a vector of $P$ observables $\mathbf y$,
\begin{equation}\label{4-6}
\mathbf y^n = \mathbf g(\mathbf u^{n}) = \begin{bmatrix}
g_1(\mathbf u^n)\\
\vdots\\
g_P(\mathbf u^n)
\end{bmatrix}, \quad 
g_p: \mathcal M\to \mathbb C \ \mbox{is an observable function $p =1,\ldots, P$}, 
\end{equation} 
Let $\mathbf g$ be restricted to an invariant subspace spanned by eigenfunctions of the Koopman operator $\mathcal K_t$. 
\end{assump}
Under this assumption, $\mathbf g$ induces a linear operator $\mathbf K$ that is finite-dimensional and advances these eigen-observable functions on this subspace~\cite{brunton2016koopman}. The physics-aware DMD Algorithm~\ref{algo:phys_dmd}~\cite{lu2019part1} can be applied to approximate the eigenvalues and eigenfunctions of $\mathbf K$ from snapshots data collected in the observable space. 

There is no principled way to choose the observables without expert knowledge of a system under consideration. Selection of observables remains a grand challenge and an active research area, e.g., machine learning and deep learning techniques were recently employed to identify the underlying Koopman operator~\cite{morton2018deep}. 
In the context of conservation laws with shocks, the equivalency between hodograph transformation and the Koopman operator, established in this study, facilitates a ``smart" choice of the observables. It is implemented via the following algorithm.

\begin{algo}{Physics-Aware DMD algorithm}\label{algo:phys_dmd}
\begin{itemize}
\item[0.] Create data matrices of $(M-1)$ observables, $\mathbf Y_1$ and $\mathbf Y_2$,
\begin{equation}\label{4-7}
\mathbf Y_1 = \begin{bmatrix}
|&|&&|\\
\mathbf y^1&\mathbf y^2&\cdots&\mathbf y^{M-1}\\
|&|&&|
\end{bmatrix}, \quad
\mathbf Y_2 = \begin{bmatrix}
|&|&&|\\
\mathbf y^2&\mathbf y^3&\cdots&\mathbf y^{M}\\
|&|&&|
\end{bmatrix}.
\end{equation}
Each column of these matrices is given by 
\begin{equation}
\mathbf y^n = \mathbf g(\mathbf u^n) = \begin{bmatrix}
g_1^n\\
g_2^n
\end{bmatrix},
\end{equation}
where $g_1^n =  x(n \Delta t, \mathbf u^n)$ is the inverse function of $u(t,x)$, and $g_2^k$ is a problem-dependent recording of shock information. 

\item[1.] Apply SVD: $\mathbf Y_1 \approx \mathbf U\boldsymbol \Sigma \mathbf V^*$, with $\mathbf U \in \mathbb C^{P\times r}, \boldsymbol\Sigma \in \mathbb C^{r\times r}, \mathbf V\in \mathbb C^{r\times M}$, and $r$ denoting the truncated rank chosen by certain criteria.
\item [2.] Compute $\tilde{\mathbf K}=\mathbf U^*\mathbf X'\mathbf V\boldsymbol\Sigma^{-1}$ as an $r\times r$ low-rank approximation for $\mathbf K$.
\item [3.] Compute eigen-decomposition of $\tilde{\mathbf K}$: $\tilde{\mathbf K} \mathbf W = \mathbf W\boldsymbol\Lambda$, $\boldsymbol\Lambda = (\lambda_k)$.
\item [4.] Reconstruct eigen-decomposition of $\mathbf K$. Eigenvalues are $\boldsymbol\Lambda$ and eigenvectors are $\boldsymbol\Phi  = \mathbf U\mathbf W$.
\item [5.] Future $\mathbf y_\mathrm{DMD}^{n+1}$ is predicted by
\begin{equation}\label{4-8}
\mathbf y_\mathrm{DMD}^{n+1} =\boldsymbol \Phi\boldsymbol\Lambda^{n+1}\mathbf b, \quad n>M%
\end{equation}
with $\mathbf b =\boldsymbol \Phi^{-1}\mathbf y_1$.
\item [6.] Transform from observables back to state-space:
\begin{equation}\label{4-9}
\mathbf u_\mathrm{DMD}^n =\mathbf g^{-1}(\mathbf y_\mathrm{DMD}^n)= g_1^{-1}(\mathbf x_\mathrm{DMD}^n).
\end{equation}
\end{itemize}
\end{algo}

\begin{remark}
Numerically, $g_1$ can be obtained by interpolation from a uniform mesh in the $(x,u)$ plane to a uniform mesh in the $(u,x)$ plane, and so can $g_1^{-1}$. The monotonicity assumption~\ref{assumption1} ensures that the observable functions are one-to-one maps.
\end{remark}

\begin{remark}
The challenge of incorporating the shock information into the Lagrangian DMD algorithm~\cite{lu2019part1} is the dependence of shock speed on the dependent variable $u$. Hodograph transformation facilitates the incorporation of this  implicitly nonlinear information by turning $u$ an independent variable and by rendering the shock speed given by the Rankine-Hugoniot condition  linear and, in fact, constant.
\end{remark}

\begin{remark}
For problems with shocks, one needs to collect snapshots until and after a shock forms. Otherwise, the Koopman operator cannot learn the shock dynamics.
\end{remark}

\begin{remark}
For mixed wave problems, one needs to collect snapshots until and after all forms of propagation occurs. This requires pre-observation, pre-processing and understanding of the data. General initial data $u_0(x)$ has to be separated into monotonic sub-regions. Physical quantities, such as shock speed and intersection point of shock and rarefaction propagation, must be understood from given data features. They give an explicit form of the shock observable function $g_2$; although problem-dependent, all the shock information is linear with respect to $u$.
\end{remark}

\section{Numerical Tests}
\label{sec:examples}

We apply the physics-aware DMD to construct ROMs of scalar conservation laws in different scenarios, including a shock, rarefaction and a mixture of both. The conservative first-order upwind scheme~\eqref{2-3} is employed as reference solution, except when an analytical solution is available. The rank truncation criterion~\eqref{2-9} with $\varepsilon = 10^{-4}$ is used in all cases.

\subsection{Riemann problem for Burgers Equation with Shock}

Consider the inviscid Burgers equation~\eqref{2-16} defined for $(x,t) \in [-0.5,1.5]\times [0,1]$ and with initial data
\begin{equation}\label{5-1}
u_0(x) = \left\{
\begin{aligned}
& 2 && \mbox{for} \ -0.5 \le x < 0\\
& 0 && \mbox{for} \ 0 \le x \le 1.5.
\end{aligned}
\right.
\end{equation}
This problem admits an analytical solution
\begin{equation}\label{5-2}
u(x,t) = \left\{
\begin{aligned}
&2&&\mbox{for}&& -0.5 \le x < st\\
&0&&\mbox{for}&& st < x \le 1.5,
\end{aligned}
\right. 
\end{equation}
where the shock speed $s = 1$ is determined from the Rankine-Hugoniot condition. The data needed for our DMD algorithm comes from the temporal snapshots of $u(x,t)$ in~\eqref{5-2}. 

The discontinuous initial data $u_0(x)$ in~\eqref{5-2} do not satisfy Assumption~\ref{assumption1}. Thus, we approximate the step function $u_0(x)$ with a smooth function, e.g., the hyperbolic tangent
\begin{equation}\label{5-3}
u_0\approx 1- \tanh\left(\frac{x}{\delta}\right), \qquad \delta \ll 1,
\end{equation}
which satisfies Assumption~\ref{assumption1}. In the $(u,x)$ plane,
\begin{equation}\label{5-4}
x_0\approx  \frac{\delta}{2}\log\left(\frac{2-u_0}{u_0}\right), \qquad \delta\ll 1.
\end{equation}
This approximation is valid in the neighborhood of the shock interface. Away from it,~\eqref{5-1} is used. Snapshots of $x(t,u)$ on a uniform mesh of $u$, which consists of $J = 2000$ equidistant points, are collected at $M = 250$ times until $T=0.25$. The ROM is used to predict the solution $u(x,t)$ for larger times, $t > T$.

\begin{figure}[tphb]
\begin{center}
\includegraphics[width=0.8\textwidth]{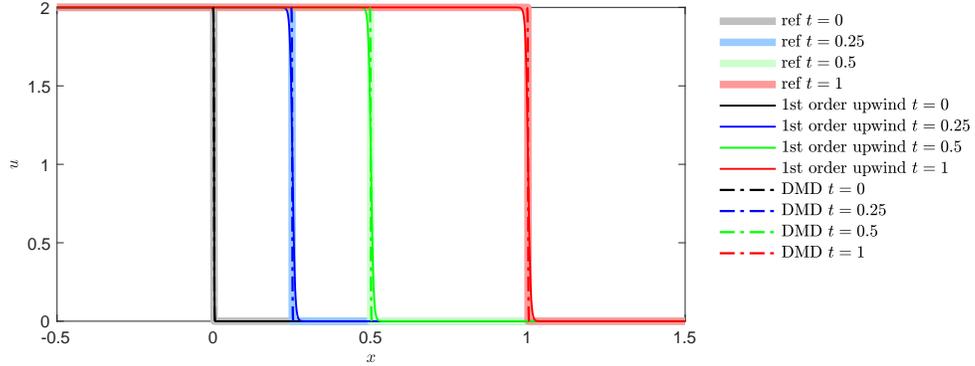}
\end{center}
\caption{Physics-aware Lagrangian DMD solution of the inviscid Burgers equation with a shock. The reference solution is given by analytic solution~\eqref{5-2}. 1st order upwind scheme by~\eqref{2-3} is also plotted here in solid line.}
\label{fig:test1}
\end{figure}

Figure~\ref{fig:test1} demonstrates the the physics-aware Lagrangian DMD algorithm with hodograph transformation captures the behavior of the shock propagation. Only $r=2$ modes are needed to construct the ROM, which remains accurate for relatively long time in the extrapolation mode. Hodograph transformation converts the nonlinear conservation law~\eqref{2-16} with discontinuous initial data~\eqref{5-1} into a linear shift with constant speed, which is readily learned from data.

\subsection{Riemann Problem for Burgers Equation with Rarefaction Wave}

Consider the inviscid Burgers equation~\eqref{2-16} defined for $(x,t) \in [-1,1]\times [0,1]$ and with initial data
\begin{equation}\label{5-5}
u_0(x) = \left\{
\begin{aligned}
& -1 && \mbox{for} \ -1 \le x < 0\\
& 1 && \mbox{for} \ 0 \le x \le 1.
\end{aligned}
\right.
\end{equation}
This problem admits an analytical solution in the form of a rarefaction wave,
\begin{equation}\label{5-6}
u(x,t) = \left\{
\begin{aligned}
&-1&&\mbox{for}&& -1 \le x < - t\\
&x/t&&\mbox{for}&&-t<x<t\\
&1&&\mbox{for}&& t < x \le 1.
\end{aligned}
\right. 
\end{equation}
A hyperbolic-tangent approximation analogous to~\eqref{5-3} is used to deal with the discontinuity in the initial data $u_0(x)$. And the same structure of data matrix is used in the DMD algorithm with $J=2000$ and $M=250$ until $t=0.25$.

\begin{figure}[tphb]
\begin{center}
\includegraphics[width=0.8\textwidth]{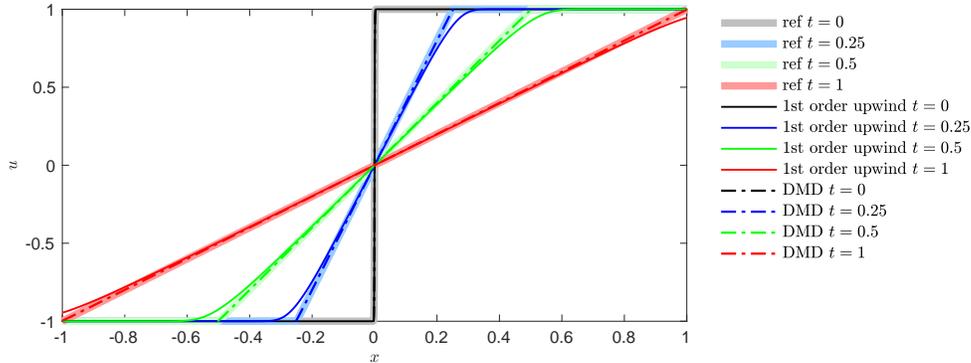}
\end{center}
\caption{Physics-aware Lagrangian DMD solution of the inviscid Burgers equation with a rarefaction wave. The reference solution is given by analytic solution~\eqref{5-6}. 1st order upwind scheme by~\eqref{2-3} is also plotted here in solid line.}
\label{fig:test2}
\end{figure}

Figure~\ref{fig:test2} shows the same satisfactory ROM results for this problem with a rarefaction wave. Only $r=2$ modes are needed to obtain  accurate predictions. As the hodograph transform, $x = x(t, u)$ satisfies a linear ODE, with $u$ acting as an independent variable. Given an accurate approximation of the initial discontinuity, i.e., selecting $\delta$ to be sufficiently small, the ROM trained on the data from generated with the analytical solution~\eqref{5-6} is even more accurate than the HFM solution. The upwind scheme~\eqref{2-3} has first-order accuracy, $O(\Delta t)$, while the Lagrangian DMD algorithm can have spectral accuracy. Visually (in Figure~\ref{fig:test1} and Figure~\ref{fig:test2}), the DMD solution has a much sharper interface than that estimated with the first-order upwind scheme~\eqref{2-3}.

\subsection{Smooth Solution of Riemann Problem for Burgers Equation with Nonmonotonic Initial Data}

Consider the inviscid Burgers equation~\eqref{2-16} defined for $(x,t) \in [0,2\pi]\times [0,1]$ and with initial data
\begin{equation}\label{5-8}
u_0(x) =1+\sin(x).
\end{equation}
Since these initial data violate Assumption~\ref{assumption1}, we decompose the interval $[0, 2\pi]$ into two parts: in the left part, $x \in [0, \pi]$, $u_0(x)$  monotonically increases; in the right part, $x \in [\pi, 2\pi]$, it monotonically decreases. Each part has a unique inverse function of $x_0(t,u_0)$. Since the shock formation time is infinite, the equation of characteristics for this Riemann problem is equivalent to~\eqref{3-5} on any finite-time interval $[0,T]$. Although this is a shock-free scenario, the two parts have different wave propagation behaviors. The numerical scheme~\eqref{2-3} with $J = 2000$ spatial discretization points and $N =1000$ time steps provides the reference solution. The data used to inform our DMD method consist of $M = 250$ snapshots of this solution.


\begin{figure}[tphb]
\begin{center}
\includegraphics[width=0.8\textwidth]{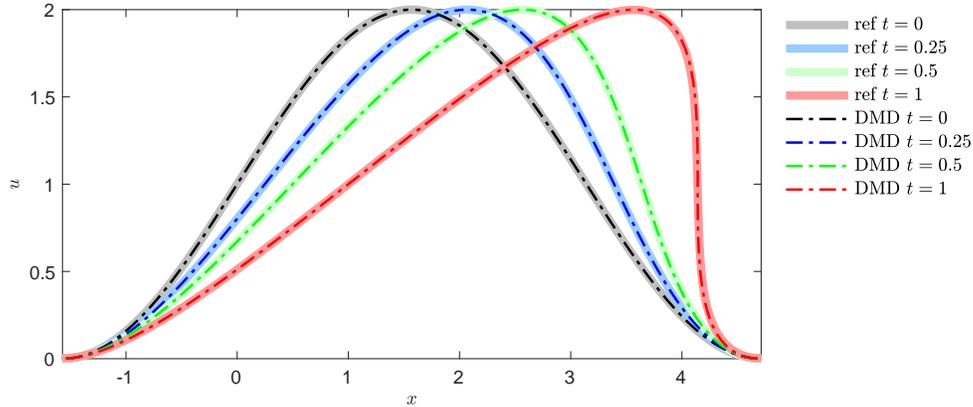}
\end{center}
\caption{Physics-aware Lagrangian DMD solution of the inviscid Burgers equation subject to nonmonotonic initial condition. The reference solution is computed with~\eqref{2-3}. }
\label{fig:test3}
\end{figure}

Figure~\ref{fig:test3} demonstrates the ability of the ROM based on our physics-aware DMD algorithm to capture these nonlinear dynamics. The ROM was trained on the early ($t \le 0.25$) data, which exhibit smooth gradients. Yet, it is capable of accurately predicting sharp gradients at later times, e.g., $t=1$. That is because, in the $(u,x)$ domain of the hodograph transform, higher gradients of $u(\cdot,x)$ translate into flatter horizontal plots of $x(\cdot, u)$.



\subsection{Riemann Problem for Burgers Equation with Rarefaction and Shock}

Consider the inviscid Burgers equation~\eqref{2-16} defined for $(x,t) \in [0,2]\times [0,1]$ and with the Gaussian-type initial data in~\eqref{2-17}. This is the setting we used to illustrate the failure of the Lagrangian DMD in section~\ref{sec:DMDfailure} (Figure~\ref{fig:goodLagrange}). The numerical scheme~\eqref{2-3} with $J =2000$ spatial discretization points and $N = 10^5$ time steps provides the reference solution. The finer time discretization is needed to satisfy the CFL constraints. The data used to inform our DMD method consist of $M = 3000$ snapshots of this solution. These data are sufficiently rich to identify the rarefaction and shock behavior of the solution. The decomposition of initial data $u_0(x)$ is needed to enforce monotonicity. The increasing branch of $u_0(x)$ is responsible for the rarefaction and its decreasing branch gives rise to the shock.

\begin{figure}[tphb]
\begin{center}
\includegraphics[width=0.8\textwidth]{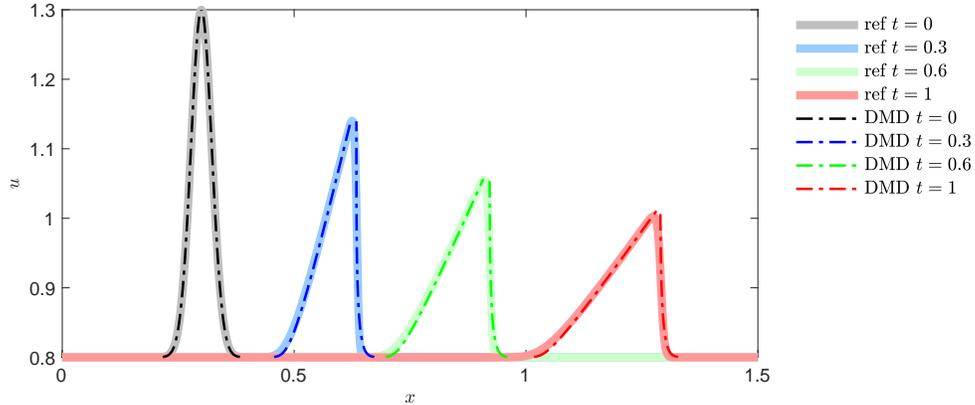}
\end{center}
\caption{Physics-aware Lagrangian DMD solution of the inviscid Burgers equation with a rarefaction wave and shock. The reference solution is computed with~\eqref{2-3}. }
\label{fig:test4}
\end{figure}

Figure~\ref{fig:test4} shows that the physics-aware DMD based on hodograph transformation provides  an accurate ROM for this Riemann problem, which could not be treated with the original physics-aware DMD. The physical shock  information, which is needed for the observable function $g_2$, includes the shock speed and the intersection point of the rarefaction wave and the shock trajectory. In this setting, the shock speed varies with time but is  still linear with respect to $u$. The physics-aware DMD algorithm can learn this linear relationship from the data with no difficulties. All of the advantages of linearity are achieved with the hodograph transform.

\subsection{Riemann problem for Buckley-Leverett Equation}

In the last numerical experiment, we consider the hyperbolic conservation law~\eqref{2-2} with a nonmonotonic flux function,
\begin{equation}
F = \frac{u^2}{u^2+a(1-u)^2}, \quad a=0.5,
\end{equation}
that is defined for $(x,t) \in [0,2]\times [0,0.5]$ and is equipped with initial data
\begin{equation}\label{5-10}
u_0(x) = \left\{
\begin{aligned}
& 1 && \mbox{for} \ 0 \leq x <1\\
& 0 && \mbox{for} \ 1 \leq x \leq 2.
\end{aligned}
\right.
\end{equation}
Similar to~\eqref{5-3}, the initial discontinuity is approximated with the hyperbolic tangent function. The hodograph treatment of this more general problem is provided in appendix~\ref{app:Buckley}.  The numerical scheme~\eqref{2-3} with $J =2000$ spatial discretization points and $N = 1000$ time steps provides the reference solution up to $t = 0.5$. The data used to inform our DMD method consist of $M = 250$ snapshots of this solution until $t = 0.125$. This set of snapshots is sufficiently rich to reveal a self-similar structure of the solution. 

Although the initial data $u_0$ are monotonic, their decomposition is needed according to the convex hull construction of the flux function (appendix~\ref{app:Buckley}). The reformulation involves two branches of different linear equations with two sets of the disjoint initial data. The shock observation function $g_2$ comprises the shock speed as well as the intersection point of the rarefaction wave and the shock trajectory. This intersection point defines the magnitude of the shock and informs the convex hull construction of the flux function. 

\begin{figure}[tphb]
\begin{center}
\includegraphics[width=0.8\textwidth]{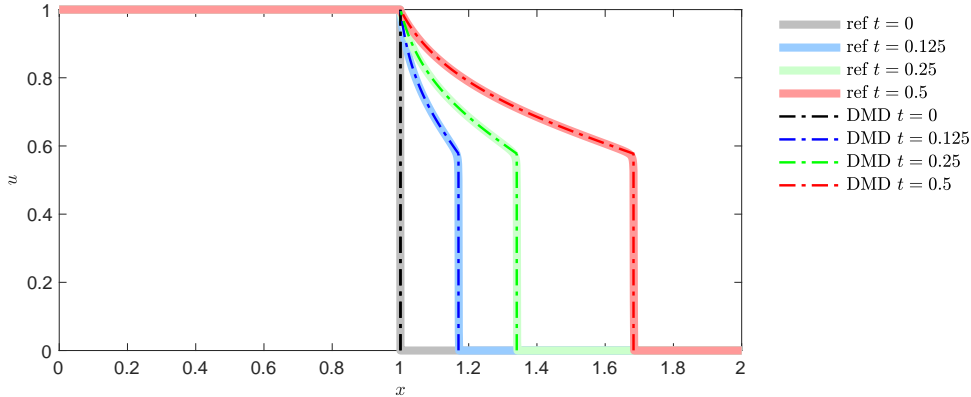}
\end{center}
\caption{Physics-aware Lagrangian DMD solution of the Buckley-Leverett equation, which has a nonmonotonic flux function. The reference solution is computed with~\eqref{2-3}. }
\label{fig:test5}
\end{figure}

Figure~\ref{fig:test5} demonstrates that the physics-aware DMD accurately captures the future states in relatively long time. Hodograph transformation allows one to determine the underlying linear Koopman operator in the nonlinear conservation laws. The iteration-free feature of DMD enhances its effectiveness and efficiency.

\section{Summary and Conclusions}
\label{sec:conc}

The Lagrangian physics-aware DMD~\cite{lu2019part1} provides a robust tool to construct ROMs of hyperbolic conservation laws, a class of problems for which standard (Eulerian) DMD methods fail. However, this algorithm is limited to problems that admit smooth strong solutions. We extended it to problems with shocks and rarefaction waves, thus addressing a long-standing challenge in ROM construction. 
This challenge stems from sever grid distortion typical of Lagrangian POD and DMD algorithms. Lacking information about shocks and discontinuities, DMD mode projection from the HFM to a ROM does not preserve the topological structure of the interface where characteristic lines cross each other. We resolved this issue by combining hodograph transformation with physics-aware DMD algorithm~\cite{lu2019part1}. 

Hodograph transforms are consistent with the Koopman operator theory in that both aim to identify linear structures in the underlying nonlinear dynamics. Our Lagrangian physics-aware DMD algorithm enhanced by hodograph transformation is capable of predicting the dynamics of weak solutions, which satisfies the entropy condition.  We demonstrated the accuracy and robustness of our algorithm on several numerical tests.

To the best of our knowledge, our study is the first to establish a connection between hodograph transformation and the Koopman operators. By providing a principled way for identifying the observables needed by the Koopman operator theory, this connection opens a door to construct ROMs for a wide range of nonlinear PDEs that are linearizable by hodograph transformation~\cite{clarkson1989hodograph}. There is an algorithmic method to do the linearization via \textit{extended} hodograph transforms. As a result, one can take advantage of the linearity and design robust iteration-free physics-aware DMD. Moreover, data-driven modeling and uncertainty quantification can be further explored using this framework. Our numerical experiments demonstrated that many physical quantities, such as the shock speed in Burgers equation and the mobility constant in Buckley-Leverett equation, can be learned from (simulation) data as long as one analyzes them in a ``smart" way.

\section{Acknowledgements}

This research was supported in part by Air Force Office of Scientific Research under award numbers FA9550-17-1-0417 and FA9550-17-1-0417, and by a gift from Total.

\begin{appendices}
\section{Scalar Conservation Laws with Convex Fluxes}
\label{app:Buckley}

Burgers' equation has a monotonically increasing flux function. Here, we extend our analysis to smooth, strictly convex flux functions $F(u)$. We consider a hyperbolic conservation law~\eqref{2-2} defined for $(x,t) \in \mathbb R \times [0,T]$. It is subject to the initial condition $u(x,0) = u_0(x)$, where the initial data $u_0(x)$ satisfy the following assumption.
\begin{assump}
The real-valued function $u_0(x)$ is such that
\begin{itemize}
\item $\lim _{x \rightarrow \pm \infty} u_0(x)=\mp 1$, and
\item $u_0(x)$ is non-increasing and, therefore, the inverse function $x(u_0)$ is well-defined on $-1 \le u_0 \le 1$.
\end{itemize}
\end{assump}

\begin{remark}
The domain of definition, $x \in \mathbb R$, can be generalized to a finite-length interval $(u_R,u_L)$. The derivation is similar.
\end{remark}

\subsection{Solution Before Shock Formation}

Similar to section~\ref{sec:3-1}, hodograph transformation yields an equation for $x(t,u)$:
\begin{equation}\label{A-2}
\frac{\text d x}{\text dt}(t, u) = f(u), \quad x(0, u) = x_0(u); \qquad u \in(-1,1).
\end{equation}
The convexity of $F(u)$ ensures that its derivative $f(u)$ is an increasing function. Let $G$ denote the inverse function of $f$:
 \begin{equation}\label{A-3}
G[f(u)] = f[G(u)] = u.
\end{equation}
Then, defining $y(t,u) = x[t,G(u)]$,~\eqref{A-2} becomes
\begin{equation}\label{A-4}
\frac{\text d y}{\text d t} (t, u) 
=u, \quad y(0,u) = y_0(u) = x_0[G(u)]; \qquad u \in(-1,1).
\end{equation}
Differentiating both sides of this equation with respect to $u$,
\begin{equation}\label{A-5}
\frac{\partial^2 y}{\partial t \partial u} = 1,
\end{equation}
which gives
\begin{equation}\label{A-6}
\frac{\text d u}{\text d t} (t, u) = y_0^{\prime}(u)+t.
\end{equation}
Therefore the shock formation time is determined by
\begin{equation}\label{A-7}
t^{*}=-\min_u y_0^{\prime}(u)=-y_0^{\prime}[f(u^{*})]
\end{equation}

\subsection{Solution After Shock Formation}

The shock speed $s$ is given by the Rankine-Hugoniot condition,
\begin{equation}\label{A-8}
s=\frac{F\left(u_{1}\right)-F\left(u_{2}\right)}{u_{1}-u_{2}},
\end{equation}
where $u_1(t)$ and $u_2(t)$ are defined as the limits of $u(t)$ from the top and bottom of the shock, respectively. Since $s = \text d x^* / \text dt$, this gives an equation for the shock trajectory $x^*(t)$,
\begin{equation}\label{A-9}
\frac{\text{d} x^*}{\text{d} t} = \frac{F\left(u_{1}\right)-F\left(u_{2}\right)}{u_{1}-u_{2}}.
\end{equation}
A system of coupled ODEs for $u_1(t)$ and $u_2(t)$ is derived in~\cite{ruiwen2018},
\begin{align}\label{A-10}
&\frac{\text{d} u_{1}}{\text{d} t} = F_1(u_1, u_2) \equiv \frac{1}{g(u_1) - f^\prime(u_1) t } \left[f(u_1) - \frac{ F(u_1) - F(u_2) }{ u_1 - u_2} \right] \\ 
&\frac{\text{d} u_{2}}{\text{d} t} = F_2(u_1, u_2) \equiv \frac{1}{g(u_2) - f^\prime(u_2) t} \left[f(u_2) - \frac{F(u_1) - F(u_2)}{u_1 - u_2}\right].
\end{align}
where $g(u) = -x'_0(u)$. These ODEs are subject to initial conditions $u_1(t^*) = u^*$ and $u_2(t^*) = u^*$.

\subsection{Summary of Hodograph Solution}
In summary, the reformulation for general scalar conservation law with convex flux is 
\begin{equation}\label{A-11}
\left\{
\begin{aligned}
&t<t^*: &&\mbox{Equation (\ref{A-2})}\\
&t>t^*:&&\left\{
\begin{aligned}
&\mbox{Equation (\ref{A-2})}&&\mbox{for}&&u\in(u_R,u_2)\cup(u_1,u_L),\\
&\mbox{Equation (\ref{A-9})}&&\mbox{for}&&u\in(u_2,u_1).
\end{aligned}
\right.
\end{aligned}
\right.
\end{equation}
where $t^* = -x'_0(u^*)$.

\begin{remark}
One can show that $u_1(t)$ is monotonically increasing in time and $u_2(t)$ is monotonically decreasing, so that
\begin{equation}\label{A-12}
u_1 \geq u^*, \quad u_2\leq u^*, \quad x'_0(u_1)+t \leq 0, \quad x'_0(u_2)+t \leq 0.
\end{equation}
In many cases of interests, and in our numerical experiment, either $u_2 = u_R$ and $u_1 = u_L$ or $|u_2 - u_1| \ll \Delta t$ (so that $u_2 \approx u_R$ and $u_1 \approx u_L$). This allows one to focus on shock propagation, i.e., on~\eqref{3-11}, without having to solve~\eqref{A-10}.
\end{remark}


\begin{remark}
For more general initial condition $u_0$, one needs to decompose $u_0(x)$ into regions of monotonicity. Each monotonic piece of $u_0$ would have a unique inverse function $x_0(u_0)$. Then, based on the generalized entropy condition, one constructs the convex hull for the flux function $F(u)$, providing a way to decompose the initial data. Shock propagating initial data and rarefaction propagating initial data are determined afterwards. Then,  the full solution is the combination of the rarefaction pieces and the shock pieces.
\end{remark}
\end{appendices}

\renewcommand\refname{Reference}
\bibliography{DMDhodograph}

\end{document}